\documentclass[11pt]{article}

\usepackage{amsmath,amsfonts,amsthm,euscript,amssymb}

\newtheorem{theorem}{Theorem}

\newtheorem{lemma}[theorem]{Lemma}

\numberwithin{equation}{section}
\numberwithin{theorem}{section}

\newcommand{\A}{\mathbb{A}}

\newcommand{\D}{\EuScript{D}}  
\newcommand{\E}{\EuScript{E}}
\newcommand{\F}{\mathbb{F}}

\newcommand{\J}{\EuScript{J}}
\newcommand{\K}{K^{\times\,2}}
\newcommand{\M}{\EuScript{M}}

\newcommand{\Q}{\mathbb{Q}}
\newcommand{\R}{\EuScript{R}}
\newcommand{\T}{\EuScript{T}}

\newcommand{\X}{\EuScript{X}}
\newcommand{\Z}{\mathbb{Z}}

\newcommand{\FF}{\EuScript{F}}
\newcommand{\KK}{\EuScript{K}}
\newcommand{\Lbig}{\Lambda_{\rm big}}

\renewcommand{\L}{\EuScript{L}}

\renewcommand{\P}{\mathbb{P}}

\renewcommand{\AA}{\EuScript{A}}

\newcommand{\End}{\mathrm{End}}
\newcommand{\Frob}{\mathrm{Frob}}
\newcommand{\GL}{\mathrm{GL}}
\newcommand{\Gal}{\mathrm{Gal}}

\newcommand{\SL}{\mathrm{SL}}
\newcommand{\Spec}{\mathrm{Spec}}

\newcommand{\ord}{\mathrm{ord}}
\newcommand{\Rad}{{\rm Rad}}
\newcommand{\Neron}{N\'eron~}
\newcommand{\etale}{\'etale~}
\newcommand{\drop}{{\rm drop}}
\newcommand{\Sp}{{\rm Sp}}
\newcommand{\tame}{\rm t}

\newcommand{\MC}{\mathrm{MC}}
\newcommand{\rk}{\mathrm{rk}}
\newcommand{\codim}{\mathrm{codim}}

\renewcommand{\div}{\mathrm{div}}

\newcommand{\Fell}{\F_\ell}
\newcommand{\Qell}{{\Q_\ell}}
\newcommand{\Zell}{{\Z_\ell}}

\newcommand{\Fqbar}{\overline{\F}_q}

\newcommand{\Kbar}{\overline{K}}

\newcommand{\kbar}{\overline{k}}
\newcommand{\taubar}{\overline{\tau}}
\newcommand{\Ubar}{U\times\kbar}
\newcommand{\Zbar}{Z\times\kbar}

\newcommand{\Cbar}{C\times \kbar}

\newcommand{\pair}[2]{\langle\,#1\!\,\,,\ \!#2\,\rangle}

\newcommand{\pr}[1]{\langle #1\rangle}
\newcommand{\gp}[1]{\langle #1 \rangle}

\newcommand{\rmk}{{\scshape Remark:~}}

\newcommand{\loccit}{{\it loc.~cit.}}

\begin{document}

\title{Big symplectic or orthogonal monodromy modulo $\ell$}
\author{Chris Hall
\thanks{Department of Mathematics, University of Texas, Austin, TX 78712, USA,
\ {\tt cjh@math.utexas.edu}}}
\maketitle

\begin{abstract}
\noindent
Let $k$ be a field not of characteristic two and $\Lambda$ be a set
consisting of almost all rational primes invertible in $k$.  Suppose
we have a variety $X/k$ and strictly compatible system $\{\M_\ell\to
X:\ell\in\Lambda\}$ of constructible $\Fell$-sheaves.  If the system
is orthogonally or symplectically self-dual, then the geometric
monodromy group of $\M_\ell$ is a subgroup of a corresponding
isometry group $\Gamma_\ell$ over $\Fell$, and we say it has big
monodromy if it contains the derived subgroup $\D\Gamma_\ell$.  We
prove a theorem which gives sufficient conditions for $\M_\ell$ to
have big monodromy.  We apply the theorem to explicit systems arising
from the middle cohomology of families of hyperelliptic curves and
elliptic surfaces to show that the monodromy is uniformly big as
we vary $\ell$ and the system.

\end{abstract}

\section{Introduction}

Galois theory in its various forms is one of the main and most powerful
tool of arithmeticians and geometers, and in particular the determination
of the Galois group of some field extensions or coverings of algebraic
varieties (or differential Galois group of some equations) is often both
an important step in solving certain problems and a very interesting
question of its own.  Whereas finding ``upper bounds'' for those Galois
groups can be relatively easy (coming from the existence of symmetries
that must be preserved), it is usually quite challenging to compute
exactly the Galois groups.  What is expected is that, given the known
symmetry constraints, the Galois group will ``usually'' be the largest
group preserving those symmetries; indeed, this is often the desired
conclusion for applications.  There are many celebrated results of this
type, among which are Serre's computations of Galois groups of torsion
fields of elliptic curves over number fields \cite{S2}, Katz's monodromy
computations in algebraic geometry leading to equidistribution statements
for angles of Kloosterman sums \cite{Ka1}, and some cases of the inverse
Galois problem \cite{Hi}, \cite{Sha1}.

It is desirable to have criteria to compute Galois groups and to show
they are ``as big as possible'', and it is most important that those
criteria involve conditions that can be checked in practice.  This paper
considers an important class of situations where the groups involved are
finite orthogonal or symplectic groups over $\Fell$.  There are quite a
few applications where such groups arise, and we describe the motivating
one in section~\ref{sec::motiv}.  Here is the simplified (and weakened)
statement of the group theoretical criterion we will prove:

\begin{theorem}
    Let $V$ be an $\Fell$-vector space together with a perfect pairing
    $V\times V\to \Fell$ and let $G\leq\GL(V)$ be an irreducible primitive
    subgroup which preserves the pairing.  If
    the pairing is symmetric,
    $G$ contains a reflection and an isotropic shear,
    and $\ell\geq 5$,
    then $G$ is one of the following:
    \begin{enumerate}
    \item the full orthogonal group $O(V)$;
    \item the kernel of the spinor norm;
    \item the kernel of the product of the spinor norm and the determinant.
    \end{enumerate}
    If
    the pairing is alternating,
    $G$ contains a transvection,
    and $\ell\geq 3$,
    then $G$ is all of the symplectic group $\Sp(V)$.
\end{theorem}

For group-theoretic terms (e.g.~transvection or isotropic shear) see
section~\ref{sec::group}.  Rather than assuming $G$ is primitive, in which
case we could appeal to \cite{Wa1} or \cite{Wa2} (cf.~section 6 of
\cite{DR}), we make explicit assumptions about a set of elements generating
$G$ and show that they (essentially) imply $G$ is primitive.  In
section~\ref{sec::group} we also give a full statement of the theorem and
its proof, and in the last two sections we give applications.  Among
those we single out (because it is easy to state and to prove) the
following (unpublished) theorem of J-K. Yu.

\begin{theorem}\label{thm0}
The mod-$\ell$ monodromy of hyperelliptic curves is $\Sp(2g,\F_{\ell})$
for $\ell>2$.
\end{theorem}

See \cite{Yu} for the preprint containing Yu's original proof or
\cite{AP} for another recent independent proof.  The theorem has
been used in several contexts.  Yu originally proved his theorem
in order to study the Cohen-Lenstra heuristics over function fields.
Chavdarov \cite{C} applied the theorem to study the irreducibility
of numerators of zeta functions of families of curves over finite
fields and Kowalski used his results to study the torsion fields
of an abelian variety over a finite field \cite{Kow2}.  Achter
applied Yu's theorem in \cite{Ac} to prove a conjecture of Friedman
and Washington on class groups of quadratic function fields.

\subsection{Acknowledgements}

We would like to thank D.~Allcock, N.M.~Katz, M.~Olsson, D.~Ulmer,
and J.F.~Voloch for several helpful conversations during the course of
research and for their interest in this work, and we would like to thank
the anonymous referree for carefully reading the paper and making several
helpful suggestions for improving the exposition.  We would also like to
thank E.~Kowalski for asking the question which motivated this paper, for
suggesting that we could extend our results, originally for orthogonal
monodromy, to symplectic monodromy, and for his generous feedback on
earlier drafts of this paper.  Finally, we must acknowledge the enormous
influence of our mentor N.M.~Katz.  This paper would not exist without
\cite{Ka3} and we hope that it will serve as a useful complement.

\subsection{Notation}

We use the notation $n\gg_{i_1,\ldots,i_m} 0$ to mean that there
is a constant $n_0(i_1,\ldots,i_m)$ depending on the objects
$i_1,\ldots,i_m$ such that $n\geq n_0(i_1,\ldots,i_m)$.  In nearly
every case the bound can be made explicit, but we usually do not
because all that matters is to determine the input objects
$i_1,\ldots,i_m$.

\section{Motivation and Strategy}\label{sec::motiv}

Let $k$ be a field, $C/k$ be a proper smooth geometrically-connected
curve, $K=k(C)$, and $E/K$ be an elliptic curve with non-constant
$j$-invariant.  In a recent paper \cite{Kow1}, E.~Kowalski asked whether
one could prove uniform big $\Fell$-monodromy results for certain
sequences of quadratic twists of $E$ when $k=\F_q$.  He went on to show
how a sufficiently strong uniform bound (as $\ell$ varies) would allow
one to prove a previously untreated variant of Goldfeld's average-rank
conjecture, although the type of bounds he required were well beyond
known results.  One piece of evidence in favor of the existence of
sufficient bounds was theorem~\ref{thm0}, an unpublished theorem due
to J.-K.~Yu.  Roughly speaking, the key property both situations share
is that the monodromy groups arise from middle cohomology of varieties.
One goal of this paper is to complete Kowalski's proof by proving bounds
of the sort he requires.  We will also show how our methods can be used
to reprove Yu's theorem.

A striking aspect of Kowalski's variant is that both $q$ \emph{and}
the degree of the conductor of the twisted curve tend to infinity in
constrast with \cite{Ka3} where only $q$ grows.  Katz fixes the degree of
the conductor for the very important reason that he wants to restrict to
nice sequences of twists which are all elements of a \emph{single} nice
geometric family.  This allows him to bring Deligne's equidistribution
theorem to bear on a variant of Goldfeld's conjecture by phrasing it in
terms of monodromy groups of $\Q_\ell$-sheaves.  Kowalski, on the other
hand, must contend with an infinite sequence of geometric families,
and he uses monodromy groups of $\Fell$-sheaves and the \v Cebotarev
density theorem as well as sieve techniques to prove his results.
A key difference between the two approaches is that, for a fixed family
of twists, the $\Qell$-monodromy groups are algebraic and essentially
independent of $\ell$ while the $\Fell$-monodromy groups are finite and
vary with $\ell$.

We fix a dense affine open $U\subset C$ and an algebraic closure
$k\to\overline{k}$.  We fix a geometric point $x\in U$, that is, an
embedding $\Spec(L)\to U$ for $L/k$ an algebraically-closed extension.  We
write $\pi_1(U)=\pi_1(U,x)$ for the \etale fundamental group and
$\pi_1^g(U)$ for the geometric fundamental group
$\pi_1(\Ubar)\leq\pi_1(U)$.  We fix a set $\Lambda$ of almost all odd
primes $\ell$ which are invertible in $k$.  For each $\ell\in\Lambda$, we
fix a lisse flat $\Z_\ell$-sheaf $\L_\ell\to U$ and let
$\rho_\ell:\pi_1(U)\to\GL_n(\Zell)$ denote the corresponding
representation.  A priori $n$ depends on $\ell$, but we assume the family
of representations $\{\rho_{\ell,\eta}=\rho_\ell\otimes\Qell\}$ is a
strictly compatible system in the sense of Serre \cite{S1}; that is, for
every $\ell\in\Lambda$, the characteristic polynomials of the Frobenii in
$\rho_{\ell,\eta}$ have coefficients in $\Q$ and are independent of $\ell$.
We write $\M_\ell\to U$ for the lisse $\Fell$-sheaf
$\L_\ell\otimes_\Zell\Fell\to U$ and say that the family $\{\M_\ell\to U\}$
is a {\it (strictly) compatible system}.

For each $\ell$, we write $G_\ell^a\leq\GL_n(\Fell)$ for the image
$(\rho_\ell\otimes\Fell)(\pi_1(U))$
and $G_\ell^g\leq G_\ell^a$ for the image of $\pi_1^g(U)$.  A priori
$G_\ell^a$ may be any subgroup of $\GL_n(\Fell)$, but if we consider
additional arithmetic information, then we may be able to deduce that
$G_\ell^a$ lies in a proper subgroup $\Gamma_\ell^a\leq\GL_n(\Fell)$.  For
example, if there is a non-degenerate pairing
$\M_\ell\times\M_\ell\to\Fell(m)$ for some Tate twist $\Fell(m)\to U$, then
we say $\M_\ell$ is {\it self dual} and we may define $\Gamma_\ell^a$ to be
the subgroup of similitudes for the pairing whose determinants are powers
of $q^m$.  One can prove a similar geometric statement: if $\M_\ell$ is
self dual and we define $\Gamma_\ell^g\leq\Gamma_\ell^a$ to be the subgroup
of isometries of the pairing, then $G_\ell^g$ lies in $\Gamma_\ell^g$.

In this paper we will assume $\M_\ell$ is self dual and
$\Gamma_\ell=\Gamma_\ell^g\leq\GL_n(\Fell)$ is the corresponding
isometry group as above for every $\ell\in\Lambda$.  We will also
assume the pairing $\M_\ell\times\M_\ell\to\Fell(m)$ is either always
symmetric or always anti-symmetric, hence $\Gamma_\ell$ is orthogonal
or symplectic respectively; recall $\ell$ is odd.  We say $\M_\ell$ has
{\it big monodromy} if $n=\mbox{rk}(\M_\ell)>1$ and $G_\ell=G_\ell^g$
contains the derived group $\D\Gamma_\ell=[\Gamma_\ell,\Gamma_\ell]$.
If $\Gamma_\ell$ is an orthogonal group, then $\D\Gamma_\ell$ is the
intersection of the kernel of the determinant and the kernel of the
spinor norm and has index four.  If $\Gamma_\ell$ is a symplectic group,
then $\Gamma_\ell=\D\Gamma_\ell$.  In particular, if we write $\Lbig$
for the $\ell\in\Lambda$ where $\M_\ell$ has big monodromy, then in either
case the index of $G_\ell$ in $\Gamma_\ell$ is uniformly bounded for
$\ell\in\Lbig$ and $\Lambda-\Lbig$ is finite, which are the sort of
properties Kowalski wants (see (16) of \cite{Kow1}).

\rmk
A priori we could relax the definition of $\Lbig$ to include all
$\ell\in\Lambda$ such that $G_\ell$ has index at most $b$ in
$\Gamma_\ell$ for some fixed $b$, but we note that the two definitions
are in fact the same for $\min\{\ell,n\}\gg_b 0$.  More precisely,
the index of the largest proper subgroup of $\D\Gamma_\ell$ grows
with $\min\{\ell,n\}$ (see table 5.2.A of \cite{KL}), so if
$\min\{\ell,n\}\gg_b 0$, then there is no proper subgroup of
$\D\Gamma_\ell$ of index at most $b$.

One of the simplest examples of a compatible system with big monodromy
can be constructed from the $\ell$-torsion of our elliptic curve
$E/K$ from above (cf.~section~\ref{geometry}).  Then $G_\ell^a$ is
the Galois group of $K(E[\ell])/K$, $G_\ell^g$ is the Galois group
of $\kbar K(E[\ell])/\kbar K$, and the function-field analogue of
Serre's theorem implies $G_\ell=\Gamma_\ell$ for almost all $\ell$;
note $\Gamma_\ell\simeq\Sp_2(\Fell)\simeq\SL_2(\Fell)$.  If we write
$g(C)$ for the genus of $C$ and $\ell\gg_{g(C)} 0$, then
$G_\ell=\Gamma_\ell$ by theorem 1.1 of \cite{CH}.  We note that
this strong uniformity was a crucial ingredient in the proof of
theorem 1.2 of that paper.

A more general example is to fix an abelian variety $A/K$ of higher
dimension with trivial $K/k$-trace and consider the compatible
system constructed from the $\ell$-torsion.  If we restrict the
endomorphism ring and dimension of $A$, then theorem 3 of \cite{S3}
implies $\ell\in\Lbig$ for $\ell\gg_A 0$, but little seems to be
known in general otherwise.  In particular, if we fix the genus of
$C$ and bound the dimensions of $A$ and its $\Kbar$-endomorphism
ring, then we do not know if $\Lbig$ may be chosen independently
of $A$.  We suspect that already for $C=\P^1$ and $\dim(A)\gg 0$
that no uniform bound exists because, roughly speaking, the
corresponding `modular varieties' are large and for arbitrarily
large $\ell$ could conceivably contain at least one line.

For general systems it is natural to ask how big $\Lbig$ is (cf.~10.7?~of
\cite{S4}).  The answer is interesting only if the Zariski closure
of $\rho_{\ell,\eta}(\pi_1^g(U))$ in $\GL_n(\Qell)$ is itself big in
an appropriate sense for any (and hence every) $\ell\in\Lambda$, so
we assume it is.  Then one can often use general methods to show that
$\Lbig$ has Dirichlet density one (see \cite{L}) or even that it contains
almost all $\ell$ (see \cite{MVW}, \cite{N} or \cite{S3}).  However,
the subset of $\ell\in\Lbig$ which these methods yield can be difficult
to describe or control, and in Kowalski's case they are insufficient
(see discussion at end of section 5 of \cite{Kow1}).  In particular,
if we let the rank $n$ tend to infinity, then these methods force us to
restrict to $\ell\gg_n 0$ where the implicit lower bound for $\ell$ tends to infinity with $n$
(e.g.~so that one can apply characteristic zero arguments).

The main goal of this paper is to demonstrate how one can prove lower
bounds for $\Lbig$ without this restriction.  The strategy we use
to achieve this is to show $G_\ell$ is an irreducible subgroup of
$\Gamma_\ell$ and then to apply theorem~\ref{thm3} where we give
sufficient criteria, in terms of a set of generators, for an irreducible
subgroup to contain $\D\Gamma_\ell$.  More precisely, we show that the
subgroup $R_\ell\leq G_\ell$ generated by the pseudoreflections is also
irreducible and use the classifications in \cite{ZS1} and \cite{ZS2}
to show that $\D\Gamma_\ell\leq R_\ell$.

For our first application we return to our last example from above
and let $\Lambda$ be the set of all odd primes $\ell$ which are
invertible in $K$.  We also let $C=\P^1$ and let $A/K$ be the
Jacobian of a curve in a special class of hyperelliptic curves
constructed in section~\ref{sec::yu}, and in theorem~\ref{thm4} we
reprove Yu's theorem and show that $\Lambda=\Lbig$.  Katz has pointed
out that the key ideas used in the proof generalize to tamely
ramified compatible systems arising from the middle cohomology of
the fibers of a Lefschetz pencil $\X\to\A^1$ of odd relative
dimension, where $\A^1=\P^1-\{\infty\}$.  Moreover, in this case
one does not need the full power of theorem~\ref{thm3} to show that
$\Lbig=\Lambda$, but instead one can appeal directly to the main
theorem of \cite{ZS2}.

The real power of theorem~\ref{thm3} emerges only when we consider more
general compatible systems.  For example, in section~\ref{sec::elliptic}
we examine systems arising from families of quadratic twists of the
elliptic curve $E/K$ when $K=\F_q(t)$.  We recall the construction due
to Katz \cite{Ka3} of an affine variety $F_d/\F_q$ which parametrizes
quadratic twists of $E/K$ by a `dense open' subset of the square-free
polynomials in $\Fqbar[t]$ of degree $d$.  We also construct, for each
$\ell$, the orthogonally self-dual lisse $\Fell$-sheaf $\T_{d,\ell}\to
F_d$ whose fibers encode the reduction modulo $\ell$ of the (unitarized)
$L$-function of the corresponding twists (cf.~section~\ref{geometry});
it corresponds to a Tate twist of the $\Qell$-sheaf constructed by Katz.
He proved that the $\Qell$-monodromy is big if $d\gg_E 0$ (cf.~theorem
1.4.3 of \loccit), and we prove something similar in theorem~\ref{thm1}:
if $\ell\geq 5$ and $d\gg_E 0$, then $\ell\in\Lbig$.

The strategy we follow to prove theorem~\ref{thm1} is to show that
the monodromy of the restriction to some one-parameter family is big.
More precisely, for each $g\in F_{d-1}$ we construct a dense open
$U_g\subset\A^1$ and a non-constant map $j_g:U_g\to F_d$ for which the
pullback $j_g^*\T_{d,\ell}\to U_g$ has big monodromy.  Up to replacing
$E/K$ by the quadratic twist $E_g/K$ and shrinking $U_g$, this reduction
amounts to restricting to the one-parameter family of twists by $c-t$
where $c\in U_g$.  We apply Katz's theory of middle convolution to analyze
the mondoromy of such a family, and in theorem~\ref{thm2} we show that
$j_g^*\T_{d,\ell}$ has big monodromy if $\ell\geq 5$ and $d\gg_{j(E)} 0$, where $j(E)$ is the $j$-invariant of $E$.  In particular, as we vary $d,g$
the collection of compatible systems $\{j_g^*\T_{d,\ell}\to U_g\}$
suffices for Kowalski's purposes.

We note that one can prove similar results for quadratic twists of
more general systems for $C$ of arbitrary genus.  Fix a dense open
$V\subset C$ and a self-dual compatible system $\{\KK_\ell\to V\}$ such
that, for each $\ell\in\Lambda$, the sheaf $\KK_\ell\to V$ is tame,
irreducible, and the monodromy around at least one geometric point of
$C-V$ is pseudoreflection.  We can construct, for each divisor $D>0$
supported on $C-V$,  a parameter space $F_D/\F_q$ of functions which is
a `dense open' subset of the Riemann-Roch space of $D$ (cf.~section 5.0
of \cite{Ka3})~and a corresponding compatible system $\{\T_{D,\ell}\to
F_D\}$.  The fibers of $\T_{D,\ell}\to F_D$ are the quadratic twists
of $\KK_\ell\to V$ and the system $\{\T_{D,\ell}\to F_D\}$ is self-dual
of the symmetry type opposite that of $\{\KK_\ell\to V\}$.  If we fix a
sufficiently nice map $t:C\to\P^1$ whose polar divisor is $[e]D$ for some
$e$, then one can argue as above to show $\T_{[de]D,\ell}\to F_{[de]D}$
has big monodromy for $\ell,d\gg_D 0$.


\section{Subgroups of Finite Symplectic and Orthogonal Groups}\label{sec::group}

Throughout this section we fix an odd prime $\ell$ and a vector space
$V$ over $\Fell$ together with a non-degenerate bilinear pairing.
We assume that the pairing is either symmetric or alternating, and in
the first case we also add the assumption that $\ell\geq 5$.  We write
$\pr{w,v}$ for the pairing of $w,v\in V$ and for a subspace $W\leq V$
we write $W^\perp$ for the orthogonal complement of $V$ and $\Rad(W)$
for the intersection $W\cap W^\perp$.  We write $\Gamma\leq\GL(V)$
for the subgroup preserving the pairing.  If the pairing is symmetric
(resp.~alternating), then we write $\Gamma=O(V)$ (resp.~$\Gamma=\Sp(V)$).

Given an element $\gamma\in \Gamma$ we write $V^{\gamma=a}$ for the
subspace of $V$ on which $\gamma$ acts as the scalar $a\in\Fell^\times$,
$V^\gamma$ for $V^{\gamma=1}$, and $V_\gamma$ for $(\gamma-1)V$.  We define
the {\it drop} of an element $\gamma\in \Gamma$ to be the codimension of
the invariant subspace $V^{\gamma=1}$.  If $\gamma\in\Gamma$ is an element
of drop 1, we say it is a {\it reflection} if $\det(\gamma)=-1$ and a
{\it transvection} if $\det(\gamma)=1$.  In either case we call $\gamma$ a
{\it pseudoreflection} and a non-zero element of $(V^{\gamma=1})^\perp$ a
{\it root}; the latter spans $V_\gamma=(V^{\gamma=1})^\perp$.  We call a
non-trivial element $\sigma\in \Gamma$ an {\it isotropic shear} if it is
unipotent and $(\sigma-1)^2=0$, and we note that the image of $\sigma-1$ is
a non-trivial isotropic subspace of $V$ and necessarily $\dim(V)\geq 4$
when the pairing on $V$ is symmetric.

\rmk If the pairing on $V$ is symmetric and $\drop(\sigma)=2$, then what
we call an isotropic shear is sometimes called a Siegel transvection.
We use the term shear in order to avoid confusion with what we call a
transvection.  It is an elementary exercise to show that there are no
(usual) transvections in the case $\Gamma= O(V)$.

\rmk What we call an isotropic shear is a quadratic element in the
sense of Thompson \cite{T}.

We say a subgroup $G\leq \Gamma$ is {\it irreducible} if $V$ is an
irreducible $G$-representation.  We say $G$ is {\it imprimitive} if $V$,
as a $G$-representation, is induced from a proper subgroup of $G$ and
otherwise we say $G$ is {\it primitive}.  We note $G$ is imprimitive if and
only if there is a non-trivial subspace $W<V$ such $V$ decomposes as a
direct sum $\oplus_{G/H}\, gW$ of the $G$-translates $gW$ over all
cosets $gH$ (cf.~section 12.D of \cite{CR}).

We devote the rest of this section to the proof of the following theorem.

\begin{theorem}\label{thm3}
    Let $r\geq 1$ and suppose $G\leq\Gamma$ is an irreducible subgroup
    together with a set of generators $S\subset G$ and a subset
    $S_0\subset S$ satisfying the following properties:
    \begin{enumerate}
    \item $\drop(\gamma)\leq r$ for every $\gamma\in S$;
    \item every $\gamma\in S-S_0$ has order prime to $(r+1)!$ {\it or is a
	pseudoreflection};
    \item $2(r+1)|S_0|\leq\dim(V)$.
    \end{enumerate}
    If the pairing is symmetric and $G$ contains a reflection and an
    isotropic shear, then $G$ is one of the following:
    \begin{enumerate}
    \item the full orthogonal group $O(V)$;
    \item the kernel of the spinor norm;
    \item the kernel of the product of the spinor norm and the determinant.
    \end{enumerate}
    If the pairing is alternating and $G$ contains a transvection,
    then $G$ is all of $\Sp(V)$.
\end{theorem}

\rmk The subgroups of $O(V)$ enumerated above are the subgroups of index
at most two excluding $SO(V)$.

If the pairing on $V$ is symmetric, then $G$ contains one or more
reflections, each of which has determinant $-1$.  Otherwise the pairing
is alternating and $G$ contains one or more transvections.  We write
$R\unlhd G$ for the normal subgroup generated by all pseudoreflections.
It is a non-trivial, although a priori it might be a proper subgroup
of $G$.  Our proof will show that it satisfies the conclusions of the
theorem, hence so does $G$.

While one can given explicit formulas for pseudoreflections in terms
of the pairing and roots, it is not necessary for what follows.  For a
fixed pseudoreflection $\gamma\in R$ most of the information about
$\gamma$ is contained in the proper subspaces $V_\gamma,V^\gamma<V$.
These spaces satisfy the key identity $V_\gamma^\perp=V^\gamma$ and
$\gp{\gamma}\leq R$ contains every pseudoreflection in $R$ with the
same (one-dimensional) root subspace.  If the pairing is symmetric,
then $\gamma$ is semisimple and $V_\gamma=V^{\gamma=-1}$, hence
$V=V_\gamma\oplus V^\gamma$.  Otherwise $\gamma$ is unipotent and it
perserves the flag $0<V_\gamma<V^\gamma<V$; the same statement is true
for an isotropic shear.

\begin{lemma}\label{lemma6}
If $W\leq V$ is a non-trivial irreducible $R$-submodule and
$H=N_G(W)\leq G$ is the stabilizer, then $\Rad(W)=0$ and
$V=\bigoplus_{G/H} gW.$
\end{lemma}

\begin{proof}
    Every $G$-translate $gW$ is an $R$-submodule because $R$ is a normal
    subgroup of $G$.  Some psuedoreflection $\gamma$ acts non-trivially on
    $W$ because otherwise $R=gRg^{-1}$ would act trivially on $gW$, hence
    on all of $V=\sum gW$, which is impossible.  Therefore the subspace of
    $W$ spanned by the roots contained in $W$ is non-trivial.  It is also
    an $R$-submodule, hence must be all of $W$, because the conjugate of a
    pseudoreflection is a pseudoreflection and so $R$ permutes the roots in
    $W$.  Thus we may write $W=\sum_\gamma W_\gamma$ where $\gamma$ varies
    over the pseudoreflections in $R$ and $W_\gamma=W\cap V_\gamma$.

    If $\gamma$ acts non-trivially on $W$, then
    $W_\gamma^\perp=V_\gamma^\perp=V^\gamma$.  Therefore, if we write
    $S\leq R$ for the subgroup generated by all pseudoreflections $\gamma$
    which act non-trivially on $W$, then $W=\sum_{\gamma\in S} W_\gamma$
    and $W^\perp=\cap_\gamma V^\gamma=V^S$.  In particular, $\Rad(W)=W\cap
    W^\perp=W\cap V^S$ is the proper, hence trivial, $R$-submodule $W^S$. 
    This proves the first part of the lemma.

    If $gW\neq W$ and $\gamma$ is a pseudoreflection, then
    $W_\gamma\cap(gW)_\gamma=0$ because $W\cap gW=0$.  Moreover, if
    $\gamma$ acts non-trivially on $W$, then $W_\gamma=V_\gamma$ and $gW$
    lies in $V^\gamma=W_\gamma^\perp$.  Therefore $gW$ lies in
    $V^S=W^\perp$ because $W=\sum_\gamma W_\gamma$, hence in general
    $g_1W\perp g_2W$ if and only if $g_1W\neq g_2W$.  In particular, the
    sum of any proper subset of $G$-translates lies in the complement of
    any unused $G$-translate, hence the sum cannot be all of $V$. 
    Therefore $V$ decomposes as the direct sum of all $G$-translates and,
    in particular, $V=\oplus_{G/H}\, gW$ because $g_1W=g_2W$ if and only if
    $g_1H=g_2H$.
\end{proof}

Our main interest in lemma~\ref{lemma6} is that it allows us to prove the
following lemma, which in turn will allow us to use classification results
about irreducible subgroups of $\GL(V)$ generated by pseudoreflections.

\begin{lemma}\label{lemma5}
$R$ is irreducible.
\end{lemma}

\rmk The argument we give below was inspired by an argument of Katz for
$\Qell$-monodromy (1.6.4 of \cite{Ka3}).

\begin{proof}
    Let $W\leq V$ and $H\leq G$ be as in the statement of
    lemma~\ref{lemma6}.  If $\dim(W)\geq r+1$, then, for every $\gamma\in
    S$, the intersection $V^\gamma\cap W$ is non-trivial because
    $\drop(\gamma)\leq r$, hence $\gamma W=W$; note, $\gamma_1 W=\gamma_2
    W$ if and only if $g_1H=g_2H$, otherwise $\gamma_1 W\cap\gamma_2 W=0$
    (cf.~(12.26) of \cite{CR}).  Therefore $W$ is stabilized by $G$ because
    $S$ generates $G$, hence $W=V$ by the irreducibility of $G$.  On the
    other hand, we cannot have $\dim(W)\leq r$ because otherwise we will
    show that it would imply $\dim(V)< 2(r+1)\cdot|S_0|$, contrary to the
    hypotheses of theorem~\ref{thm3}.  To prove this we need two
    lemmas.

\begin{lemma}\label{lemma7}
    For every $\gamma\in S$, if $\{g_iW\}$ is a subset of at least
    $(r+1)/\dim(W)$ $G$-translates, then $\{g_iW\}$ contains a
    $\gp{\gamma}$-orbit.
\end{lemma}

\begin{proof}
	The subspace $\bigoplus_i g_iW$ is at least $(r+1)$-dimensional,
	so it intersects $V^\gamma$ non-trivially.  Suppose $v\neq 0$
	lies in the intersection.  The non-empty subset of translates
	in $\{g_iW\}$ such that the projection of $v$ onto $g_iW$
	is non-trivial is $\gp{\gamma}$-stable, hence it is a union
	of $\gp{\gamma}$-orbits.
\end{proof}

    For any $\gamma\in S$ we say a $\gp{\gamma}$-orbit in $\{gW\}$
    is non-trivial if it has at least two elements.  If $\gamma$
    is a pseudoreflection, then every $\gp{\gamma}$-orbit is trivial
    because every $G$-translate $gW$ is an $R$-module and $\gamma\in R$.
    If $\gamma\in S-S_0$ is not a pseudoreflection, then a non-trivial
    $\gp{\gamma}$-orbit would have to contain at least $r+2$ elements,
    which is impossible by lemma~\ref{lemma7}.  Therefore $\gp{\gamma}$
    acts trivially on $\{gW\}$ for every $\gamma\in S-S_0$.

\begin{lemma}\label{lemma8}
    Suppose $\dim(W)\leq r$ and $\gamma\in S_0$.  If the $i$th
    $\gp{\gamma}$-orbit in $\{gW\}$ has $e_i$ elements, then
    $\dim(W)\sum_i (e_i - 1)< r+1$.
\end{lemma}

\begin{proof}
	This follows immediately by considering any subset of $\{gW\}$
	containing at most $e_i-1$ elements from the $i$th orbit.
	In particular, such a subset contains no $\gp{\gamma}$-orbit,
	hence lemma~\ref{lemma7} implies it has less than $(r+1)/\dim(W)$
	elements.
\end{proof}

    By applying lemma~\ref{lemma8} to $\gamma\in S_0$ and the set of
    non-trivial $\gp{\gamma}$-orbits in $\{gW\}$ we conclude that there
    are less than $(r+1)/\dim(W)$ non-trivial $\gp{\gamma}$-orbits
    and they have less than $2(r+1)/\dim(W)$ elements in total.
    In particular, the union of all such orbits for all $\gamma\in
    S_0$ has less than $2(r+1)|S_0|/\dim(W)$ elements.  On the other
    hand, $S$ generates $G$ and $G$ is irreducible, so every conjugate
    $gW$ lies in a non-trivial orbit for some $\gamma\in S_0$, hence
    $\dim(V)<2(r+1)|S_0|$.  Therefore $W=V$.
\end{proof}

If the pairing on $V$ is alternating, then lemma~\ref{lemma6} together
with the main theorem of \cite{ZS2} imply that $R=\Sp(V)$; for $\dim(V)=2$
this is a well-known result of \cite{D}.  Therefore we may assume $R$
is generated by reflections and apply the classification of irreducible
reflection groups in \cite{ZS1}.  There are a handful of exceptional
groups, but the key to eliminating them is to find a pair of reflections
such that the order of their product is sufficiently large.

\begin{lemma}\label{lemma9}
There are conjugate reflections $\rho,\sigma\rho\sigma^{-1}\!\in R$ such
that $[\rho,\sigma]=\rho\sigma\rho\sigma^{-1}$ has order at least $\ell$.  
\end{lemma}

\begin{proof}
    First, we claim that for every isotropic shear $\sigma\in G$ there
    is a reflection $\rho\in R$ such that $\rho\sigma\neq\sigma\rho$.
    For the roots of all reflections in $R$ span $V$ by the
    irreducibility of $R$ and we claim $\rho\sigma=\sigma\rho$ if and
    only if $V_\rho\leq V^\sigma$.  In particular, at least one root does
    not lie in $V^\sigma$ and so the corresponding reflection satisfies
    $\rho\sigma\neq\sigma\rho$ as claimed.  To prove the claim we observe
    that $V_\rho\leq V^\sigma$ and $V_\sigma\leq V^\rho$ are dual, with
    respect to the pairing, hence both statements hold or neither does.
    The first implies that $(\sigma-1)(\rho-1)=0$ in $\End(V)$ and the
    second that $(\rho-1)(\sigma-1)=0$.  If both hold, then clearly
    $\sigma\rho=\rho\sigma$.  Conversely, if $\rho\sigma=\sigma\rho$,
    then $\rho\sigma(r)=-\sigma(r)$ for any root $r\in V_\rho$, hence
    $\sigma$ stabilizes $V_\rho$ and its restriction acts trivially,
    so $V_\rho\leq V^\sigma$.

    Next, fix any pair $\rho,\sigma$ which do not commute, a root $z$ of
    $\rho$, and $y\in V^\sigma\backslash V^H$.  Note, $z\not\in V^\sigma$
    because $V_\rho\not\leq V^\sigma$, so $y,z$ are independent.  Rescaling
    $y$ if necessary we may assume $\rho y=y+z$.  We claim that
    $x=\sigma z-z$ lies in $V^H$ and hence $x,y,z$ span a three-dimensional
    $H$-submodule $W$.  By definition $x$ lies in the isotropic subspace
    $V_\sigma\leq V^\sigma$, hence $\pair{x}{x}=0$, and to show that it
    lies in $V^\rho$ it suffices to show $\pair{x}{z}=0$.  Using the
    identities $\pair{x}{x}=0$ and $\pair{\sigma z}{\sigma z}=\pair{z}{z}$
    one easily deduces that $\pair{\sigma z}{z}=\pair{z}{z}$, hence
    $\pair{x}{z}=\pair{\sigma z}{z}-\pair{z}{z}=0$ as claimed.

    Finally, in terms of the ordered basis $x,y,z$ of $W$ we have
    \[
        \rho = \left(\begin{array}{rrr}
            1 & 0 & 0 \\ 0 & 1 & 1 \\ \,0 & \,0 & \!\!-1
        \end{array}\,\right)\!,\quad 
        \sigma = \left(\begin{array}{rrr}
            1 & 0 & 0 \\ 0 & 1 & 0 \\ 1 & 0 & 1
        \end{array}\,\right)\!,\quad 
        \rho\sigma\rho\sigma^{-1} = \left(\begin{array}{rrr}
            1 & 0 & 0 \\ 1 & 1 & 0 \\ \!\!-2 & 0 & 1
        \end{array}\,\right) \in\GL(W).
    \]
    Thus the restriction $u$ of $[\rho,\sigma]$ to $W$ satisfies
    $(u-1)^\ell=0$, hence $u$ has order $\ell$ and the lemma follows.
\end{proof}

Applying lemma~\ref{lemma9} to the classification in \cite{ZS1} we find
either $[O(V):R]\leq 2$ or $R$ is imprimitive.  In their notation the
two families of imprimitive groups we must rule out are $G(m,m,n)$ and
$G(2m,m,n)$, where $n=\dim(V)$.  The first is the subgroup of $\GL(V)$
generated by the permutation matrices and the diagonal matrix with
diagonal $(\zeta,1/\zeta,\ldots,1)$, where $\zeta$ is a primitive
$m$th root of unity in $\Fell^\times$ and $m>1$.  The second is the
group generated by $G(2m,2m,n)$ and the diagonal matrix with diagonal
$(-1,1,\ldots,1)$.  We can eliminate all the imprimitive groups but
$G(2,2,n)$ and $G(2,1,n)$ using the following lemma.

\begin{lemma}
If $G(m,m,n)\leq O(V)$, then $m=2$.  If $G(2m,m,n)\leq O(V)$, then $m=1$.
\end{lemma}

\begin{proof}
    $G(2m,2m,n)$ is a subgroup of $G(2m,m,n)$, hence the first
    statement of the lemma implies the second.  Let $\pi$ be any
    permutation matrix in $O(V)$ relative to some basis $\{x_i\}$
    of $V$.  Then $\pair{\pi x_i}{\pi x_j}=\pair{x_i}{x_j}$,
    hence $\pair{x_i}{x_j}=a\delta_{ij}+b(1-\delta_{ij})$ for some
    constants $a,b\in\Fell$, where $\delta_{ij}$ is the Kronecker delta
    function.  Let $M\in O(V)$ be a diagonal matrix whose diagonal is
    $(\zeta,1/\zeta,1,\ldots,1)$, where $\zeta$ is a primitive $m$th root
    of unity.  We see that $b=0$ because $b=\pair{M x_1}{M x_3}=\zeta
    \pair{x_1}{x_3}=\zeta b$ and $\zeta^2=1$ because $a=\pair{ M x_1}{M
    x_1}=\zeta^2 a$.
\end{proof}

To eliminate the two remaining imprimitive groups we note that the product
of any pair of rotations in either of these groups has order at most 4,
hence we can eliminate them using lemma~\ref{lemma9}.  Therefore the only
possibilities for $R$ remaining are the groups of index at most two in
$O(V)$ which are not $SO(V)$, which completes the proof of the theorem.

\section{A Theorem of Yu}\label{sec::yu}

Let $q$ be an odd prime power, $C=\P^1$ over $\F_q$ and $K$ be the global
field $\F_q(C)=\F_q(t)$.  Fix $g\geq 1$ and a monic square-free $f(x)\in
\F_q[x]$ of degree $2g$.

Let $X/K$ be the hyperelliptic curve which is the natural (one-point)
compactification of the affine curve $y^2=(t-x)f(x)$.  The Jacobian
$J/K$ of $X$ is a $g$-dimensional abelian variety and for any rational
prime $\ell$ not dividing $q$ we write $J[\ell]$ for the subgroup of
$\ell$-torsion.  The main goal of this section is to prove the following
theorem due to Jiu-Kang Yu \cite{Yu}.

\begin{theorem}\label{thm4}
    If $\ell$ is odd, then the group $G_\ell=\Gal(K(J[\ell])/K)$ is as big
    as possible.  More precisely, there is a primitive $\ell$th root of
    unity $\zeta_\ell\in K(J[\ell])$ and $K(J[\ell])/K(\zeta_\ell)$ is a
    geometric extension with Galois group $\Gamma_\ell=\Sp(2g,\Fell)$.
\end{theorem}

\rmk For $g=1$ the theorem is equivalent to Igusa's theorem \cite{I}
for the so-called Legendre curve.

Our proof, which will occupy the remainder of this section, differs from
the proof in \cite{Yu} and has the advantage that the techniques used
allow one to prove more general results.

By the existence of the Weil pairing we know that $\mu_\ell\subset
K(J[\ell])$, so to prove the theorem we may make the finite scalar
extension where we replace $K$ by $K(\mu_\ell)$.  If we fix an isomorphism
$\mu_\ell\simeq\Fell$, then the group law together with the Weil pairing
gives $J[\ell]$ the structure of a $2g$-dimensional $\Fell$-vector
space together with a non-degenerate alternating pairing.  Therefore if
we choose a basis of $J[\ell]$, then we may identify $G_\ell$ with a
subgroup of $\Gamma_\ell$ and we must show that $G_\ell=\Gamma_\ell$.

Let $\X\to C$ denote the minimal regular model of $X/K$ and $\J\to C$
the \Neron model of $J/K$.  The fibers of $\X\to C$ are proper smooth
curves of genus $g$ over the open complement $j:U\to C$ of the finite
subset $Z=\{\tau\in\Fqbar:f(\tau)=0\}\cup\{\infty\}$, and each fiber
of $\X\to Z-\{\infty\}$ is smooth away from an ordinary double point
(i.e.~Lefschetz).  In particular, the restriction of $\J\to C$ to
$\A^1=C-\{\infty\}$ has semistable reduction.  Over $\infty$ the fibers
of $\X,\J$ are more difficult to describe, but we will show they are
sufficiently `tame' and that we can ignore them.

Multiplication by $\ell$ on $J$ extends to an isogeny of $C$-group
schemes $\times\ell:\J\to\J$ and we define $\J_\ell\subset\J$
to be the kernel.  The latter is a quasi-finite \etale group
scheme over $C$ and the restriction $\J_\ell\to U$ is finite
\etale\!\!.  If we write $\pi:\X\to U$ for the restriction of
$\X\to C$, then (the sheaf of sections of) $\J_\ell\to C$ is
isomorphic to the direct image sheaf $j_*\R^1\pi_*\mu_\ell$.
More precisely, the fiber of $\R^1\pi_*\mu_\ell\to U$ over a
geometric point $\taubar\in U$ is the \etale cohomology group
$H^1(C_{\taubar},\mu_\ell)=\J_{\ell}(\F_q(\taubar))$ and the adjunction
map $j_*\R^1\pi_*\mu_\ell\simeq j_*j^*\J_\ell\to \J_\ell$ is an
isomorphism.

We fix a geometric generic point $\taubar\in U$ (i.e.~an algebraic
closure $\Kbar/K$) and let $\pi_1(U)=\pi_1(U,\taubar)$ denote
the \etale fundamental group.  The lisse sheaf $\J_\ell\to U$
corresponds to a $\Fell$-representation of $\pi_1(U)$ on the fiber
$(\J_\ell)_{\taubar}=J[\ell]$ (which is defined up to inner automorphism).
More precisely, if we write $G_K=\Gal(\Kbar/K)$, then the quotient $G_K\to
G_\ell$ factors through $G_K\to\pi_1(U)$.  In fact, the following lemma
implies $\pi_1(U)\to G_\ell$ factors through the maximal tame quotient
$\pi_1(U)\to\pi_1^{\tame}(U)$.

\begin{lemma}\label{lemma1} 
$\J_\ell\to C$ is tamely ramified.
\end{lemma}

\begin{proof}
    The extension $K(J[2])/K$ is a scalar extension hence unramified
    and by Kummer theory $K(J[4])/K(J[2])$ is tamely ramified, hence
    $K(J[4])/K$ is tamely ramified.  Raynaud's criterion for semi-stable
    reduction implies $J$ has semi-stable reduction over $L=K(J[4])$
    (cf.~4.7 of \cite{G2}), hence $L(J[\ell])/L$ is tamely ramified and
    thus so is $K(J[\ell])/K$.
\end{proof}

We fix an algebraic closure $\F_q\to\Fqbar$ and let
$\pi_1^{\tame}(\Ubar)\leq\pi_1^{\tame}(U)$ denote the geometric subgroup. 
If we order the points in $\Zbar$, then for each $c\in\Zbar$ we may choose
a topological generator $\sigma_c$ of the inertia group
$I(c)\leq\pi_1^{\tame}(\Ubar)$ so that the ordered product is the identity
(cf.~\cite{G1} or \cite{SGA1}).  Moreover, $\pi_1^{\tame}(\Ubar)$ is
topologically generated by $\sigma_c$ for $c\in \Zbar-\{\infty\}$, hence to
prove the theorem it suffices to show that the images of these elements
generate $\Gamma_\ell$ which we will do using theorem~\ref{thm3}.  We note
that it follows $K(J[\ell])/K$ is geometric (i.e.~$\kbar\cap K(J[\ell])=k$)
because the image of $\pi_1^{\tame}(U)$ in $\Gamma_\ell$ lies between the
image of $\pi_1^{\tame}(\Ubar)$ and  $\Gamma_\ell$, hence is $\Gamma_\ell$.

The Picard-Lefschetz formulas imply that $\sigma_c$ acts as a
(symplectic) transvection on $J[\ell]$ for every $c\in\Zbar-\{\infty\}$
(cf.~theorem III.4.3 of \cite{FK}).  One can also use Katz's theory
of middle convolution \cite{Ka2} for $\Fell$-sheaves to deduce the
same thing as well as to give another proof of lemma~\ref{lemma1}.
More importantly, one can also show that $\J_\ell\to C$ is irreducible
(i.e.~$\pi_1^{\tame}(\Ubar)$ acts irreducibly on $J[\ell]$).
The key is to identify $\J_\ell\to C$ with the middle convolution
$\MC_{-1}(\L_{\chi(f(x))})$ which is irreducible (see next section for
notation and details).

To complete the proof of theorem~\ref{thm4} we let $\Gamma=\Gamma_\ell$,
$G=G_\ell$, $r=1$, $S=\{\sigma_c:c\in \Zbar-\{\infty\}\}$, and
$S_0=\emptyset$ and apply theorem~\ref{thm3}.   Note the image
of $\sigma_c$ in $G_\ell$ has prime order $\ell>2$ for every
$c\in\Zbar-\{\infty\}$.

\subsection{Middle Convolution}\label{sec::conv}

Let $C=\P^1$, $K=k(x)$, and assume $k=\kbar$.  Let $\FF\to C$ be
a tame quasi-finite \'etale sheaf with coefficients in $\F_\ell$.
We say $\FF$ is N\'eron if there is a dense open set $j:U\to C$
so that the restriction $\FF\to U$ is lisse and the adjunction map
$j_*j^*\FF\to\FF$ is an isomorphism.  Let $\M$ denote the collection
of irreducible N\'eron $\FF\to C$ such that the generic rank of $\FF$
is at least two or the restriction $\FF\to C-\{\infty\}$ has at least
two ramified fibers.  If $U\subset C-\{0,\infty\}$ is dense open, then
let $\M_U$ denote the subset of $\FF\in\M$ such that the restriction
$\FF\to U$ is lisse and let $\rk(\FF)$ denote the rank of $\FF\to U$.
Moreover, if $t\in C-U$, then let $\FF(t)$ denote the representation of
$I(t)\leq\pi_1^t(U)$ corresponding to $\FF\to U$.  We fix an ordering
of $C-U$ and topological generators $\sigma_t\in I(t)$ for $t\in C-U$
so that the ordered product is the identity in $\pi_1^t(U)$.

If $\lambda\in\F_\ell^\times$ has order invertible in $k$ and $t\in C-U$,
then we write $\lambda:I(t)\to\F_\ell^\times$ for the representation
$\sigma_t\mapsto\lambda$.  Let $\L_\lambda=\L_{\chi(x)}\to C$ denote
the Kummer sheaf whose restriction to $C-\{0,\infty\}$ is lisse and
invertible and for which $\L_\lambda(0)=\lambda$; note, $\L_\lambda$
is tame and N\'eron but does not lie in $\M$.  If $t\in C-\{\infty\}$
is a closed point and $\tau_t:C\to C$ is the involution $x\mapsto t-x$,
then $\tau_t^*\L_\lambda$ is lisse over $C-\{t,\infty\}$.  Moreover,
if $\FF\in\M_U$ and $i$ denotes the open immersion $U-\{t\}\to C$,
then sheaf $i_*i^*(\FF\otimes\tau_t^*\L_\lambda)$, the twist of $\FF$
by $\tau_t^*\L_\lambda$, also lies in $\M$.  In particular, the \etale
cohomology group $V_t=H^1(C,i_*i^*(\FF\otimes\tau_t^*\L_\lambda))$ is
the only non-vanishing group and its dimension is constant as $t\in U$
varies over the closed points.

The {\it middle convolution} of $\FF\in\M_U$ by $\L_\lambda$, which we
denote $\MC_\lambda(\FF)$, is the sheaf in $\M_U$ whose fiber over $t$
is $V_t$.  Katz defined it when the characteristic of $k$ is positive
(cf.~\cite{Ka2}), and one can use \cite{DR} to extend the definition to
characteristic zero (where everything is tame).

\begin{lemma}\label{lemma12}
The `functors' $\MC_\lambda:\M_U\to\M_U$ satisfy:
\begin{enumerate}
\item $\MC_1(\FF)\simeq\FF$;
\item $\MC_{\lambda_1}(\MC_{\lambda_2}(\FF))
	    \simeq\MC_{\lambda_1\lambda_2}(\FF)$;
\item $\rk(\MC_\lambda(\FF))
	    = \sum_{t\in\A^1-U}\codim\left(\FF(t)^{I(t)}\right)
		- \dim\left((\FF(\infty)\otimes\lambda)^{I(\infty)}\right)$;
\item $\MC_\lambda(\FF)(t)/\MC_\lambda(\FF)(t)^{I(t)}
	    \simeq \left(\FF(t)/\FF(t)^{I(t)}\right)\otimes \lambda$
    for $t\in\A^1-U$.
\end{enumerate}
\end{lemma}

\begin{proof}
    The first two statements show that $\MC_\lambda$ is `multiplicative'
    in $\lambda$.  See proposition 2.9.7 of \cite{Ka2} or proposition
    3.2 and theorem 3.5 of \cite{DR}.  Statements 3 and 4 follow from
    corollary 3.3.6 of \cite{Ka2} or from lemma 2.7 and lemma 4.1
    respectively of \cite{DR}.
\end{proof}

If $t\in\A^1-U$ and $F=\FF(t)$, then the Jordan decomposition of $F$
together with the above properties completely determines the decomposition
of $M=\MC_\lambda(\FF)(t)$.  Let $U_n$ denote an irreducible unipotent
Jordan block of size $n$.  The number of unipotent (or trivial) blocks of
the form $U_1$ in either $F,M$ is the dimension of the respective space
of $I(t)$-invariants, and there is a bijection between the non-trivial
blocks of $F$ and those of $M$:
\begin{equation}\label{eqn1}
    U_n\mapsto U_{n-1}\otimes\lambda,
    \quad U_n\otimes 1/\lambda\mapsto U_{n+1},
    \quad B\not\in\{ U_n,U_n\otimes 1/\lambda\} \mapsto
		B\otimes\lambda \not\in\{U_n\otimes \lambda,U_n\}.
\end{equation}
The dimension changes in the first two cases are due to the isomorphism
$U_m\simeq U_{m+1}/U_{m+1}^{I(t)}$ applied to unipotent blocks of $F,M$
respectively.

For example, let $g\geq 1$, $f(x)\in k[x]$ be square-free of degree
$2g$, and  $\FF\in\M$ be the quadratic Kummer sheaf $\L_{\chi(f(x))}$.
If we write $U=\A^1-\deg_0(f)$, then $\FF\in\M_U$ and the Tate twist
$\MC_{-1}(\FF)(1)$ is the N\'eron sheaf $\J_\ell\to C$ of the previous
section.  More precisely, if $t\in U$ is  a closed point, then the fiber
of $\MC_{-1}(\FF)(-1)$ over $t$ is $H^1(C,\L_{\chi(f(x)(t-x))}(1))$, and
the latter is easily seen to be cohomology group $H^1(X_t,(\Z/\ell)(1))$
of the hyperelliptic curve $X_t/k$.  For each $t\in\A^1-U$, the
quotient $\FF(t)/\FF(t)^{I(t)}$ is the scalar representation $-1$, so
$\MC_{-1}(\FF)(t)/\MC_{-1}(\FF)(t)^{I(t)}$ is the trivial representation
$\Z/\ell$.  Thus, $\MC_{-1}(\FF)(t)$ has one Jordan block of the form
$U_2$ and the rest are all trivial, hence the monodromy is a transvection.
From formula 3 of lemma~\ref{lemma12} we see that $\rk(\MC_{-1}(\FF))=2g$;
note, $I(\infty)$ acts trivially on $\FF(\infty)$.

\section{Quadratic Twists of Elliptic Curves}\label{sec::elliptic}

Let $q$ be an prime power not divisible by 2,3 and fix an algebraic
closure $\F_q\to\Fqbar$, although we remark that most of the results
in this section apply if we replace $\F_q$ by an arbitrary field of
characteristic distinct from 2,3.  We fix a proper smooth geometrically
connected curve $C/\F_q$ and write $K=\F_q(C)$ for its function field.

\subsection{Geometry of a Twisted Curve}\label{geometry}

We fix an elliptic curve $E_1/K$ with non-constant $j$-invariant and
write $\E_1\to C$ for its \Neron model.  For every non-trivial coset
$f\K\subset K^\times$ we write $E_f/K$ for the so-called quadratic twist
of $E_1/K$ by $f$.  It is the unique elliptic curve over $K$ which is not
$K$-isomorphic to $E_1$ but is $K(\sqrt{f})$-isomorphic.  The \Neron model
$\E_f\to C$ of $E_f/K$ is a smooth group scheme and the group of sections
$\E_f(C)$ is canonically isomorphic to the Mordell-Weil group $E_f(K)$.

Let $\ell$ be a prime which is invertible in $K$.  The multiplication by
$\ell$ map on $E_f(K)$ extends to an isogeny $\times\ell:\E_f\to\E_f$ and
we define $\E_{f,\ell}\subset\E_f$ to be the kernel.  It is a quasi-finite
\etale group scheme over $C$ which we call the \Neron model of the
$\ell$-torsion $E_f[\ell]$.  We note that if $j:U\to C$ is the inclusion
of a non-empty open set, then $\E_{f,\ell}$ is canonically isomorphic to
$j_*j^*\E_{f,\ell}$, hence is a so-called middle extension.  The fiber
of $\E_{f,\ell}$ over a geometric generic point of $C$ is $E_f[\ell]$ and
over an arbitrary geometric point of $C$ it is a subspace of $E_f[\ell]$.

\bigskip
\begin{lemma}
    Let $v$ be a geometric point of $C$.  If $\ell$ does not divide
    the order of the component group of the special fiber of $\E_f$
    over $v$, then
    $$  \dim(\E_{f,\ell}(v)) = 
        \begin{cases}
           2 & \mbox{if $\E_f$ has good reduction over $v$} \\
           1 & \mbox{if $\E_f$ has multiplicative reduction over $v$} \\
           0 & \mbox{if $\E_f$ has additive reduction over $v$} \\
        \end{cases}.
    $$
\end{lemma}

\begin{proof}
    The assumption on the order of the component group over $v$
    ensures that all $\ell$-torsion lies in the identity component
    of the special fiber.  The lemma follows easily by observing
    that identity component is cyclic in the case of multiplicative
    reduction and $\ell$-torsion free in the case of additive reduction
    (cf.~proposition 5.1 of \cite{Si}).
\end{proof}

For $\ell$ satisfying the hypothesis of the lemma we regard $\E_{f,\ell}$
as the $\Fell$-analogue of Kodaira's homological invariant (cf.~section
7 of ~\cite{Kod} and section 2 of \cite{Shi}).  In order to obtain a
similar result for general $\ell$ one must restrict to the intersection
of $\E_{f,\ell}$ with the identity component of $\E_f$.  However, for
ease of exposition we will assume $\ell\neq 2,3$ and $\ell$ does not
divide $\max\{1,-\ord_v(j(E_f))\}$ for every closed point $v\in C$,
where $j(E_f)=j(E_1)$ is the $j$-invariant, so that $\E_{f,\ell}$ is
`connected'.  We note that the set of exceptional $\ell$ is independent
of $f$.

We say that $\E_{f,\ell}$ has big monodromy if the Galois group of
$K(E_f[\,\ell\,]\,)/K$ contains $\SL_2(\Fell)$.  For $\ell\geq 5$
we note that there are no index-two subgroups of $\SL_2(\Fell)$, hence
$K(E_f[\ell])$ and $K(\sqrt{f})$ are (geometrically) disjoint extensions,
so $\E_{f,\ell}$ has big monodromy if and only if $\E_{1,\ell}$ does.
In particular, after we expand the set of exceptional $\ell$ to include
those such that $\E_{f,\ell}$ does \emph{not} have big monodromy we obtain
a finite set which is still independent of $f$, and we assume $\Lambda$
is a subset of the complement.

\rmk If we write $g(C)$ for the genus of $C$, then theorem 1.1 of
\cite{CH} asserts $\E_{f,\ell}$ has big monodromy if $\ell\gg_{g(C)}0$,
hence the set of exceptional $\ell$ may even be taken to depend
only on $g(C)$ and the primes dividing the coefficients of the
divisor of poles of $j(E_1)$.

We write $M_f,A_f\subset C$ for the divisors of multiplicative
and additive reduction respectively of $\E_{f,\ell}\to C$ and let
$U_f=C-M_f-A_f$.  By assumption $\E_{f,\ell}$ has big monodromy, and in
particular, the restriction of $\E_{f,\ell}$ to $U_f$ is an irreducible
lisse $\Fell$-sheaf of rank two.

\begin{lemma}\label{lemma11}
    The \etale cohomology groups of $\E_{f,\ell}$ over $C\times\Fqbar$ are
    $\Fell$-vector spaces satisfying
    $$
	\dim(H^i(C\times\Fqbar,\E_{f,\ell})) =
	    \begin{cases}
		\deg(M_f)+2\cdot\deg(A_f)+4\cdot({\rm genus}(C)-1)
		    & \mbox{if $i=1$} \\
		0 & \mbox{otherwise}
	    \end{cases}.
    $$
\end{lemma}

\begin{proof}
    The cohomology groups for $i\neq 0,1,2$ are trivial by the
    cohomological dimension of $C$.  The groups are also trivial
    for $i=0,2$ because $\E_{f,\ell}$ is irreducible of rank two:
    these groups are the $\pi_1(U_f\times\Fqbar)$-invariants and
    $\pi_1(U_f\times\Fqbar)$-coinvariants respectively of $E_f[\ell]$
    (because $\E_{f,\ell}$ is a middle extension and is irreducible of
    rank greater than one).  Finally, the dimension for $i=1$ follows from
    the (N\'eron--)Ogg--Shafarevich formula (see \cite{O} or \cite{Sha2});
    note, $\E_{f,\ell}$ is tamely ramified because 6 is invertible in $K$.
\end{proof}

Let $V_{f,\ell}$ be the \etale cohomology group
$H^1(C\times\Fqbar,\E_{f,\ell})$.

$\Frob_q$ acts on $V_{f,\ell}$ by functoriality, hence we may regard
$\Frob_q$ as a conjugacy class of elements in $\GL(V_{f,\ell})$.  On the
other hand, if $L(T,E_f/K)$ is the $L$-function of $E_f/K$, which we
note lies in $\Z[T]$ and has leading coefficient which is non-zero modulo
$\ell$, then
    $$ \det(1-qT\Frob_q|H^1(C\times\Fqbar,\E_{f,\ell}))
	    \equiv L(T,E_f/K)\pmod{\ell}.
    $$
In particular, we obtain an upper bound on the order of vanishing of
$L(T,E_f/K)$ at $T=1/q$, the so-called analytic rank of $E_f/K$, by
studying the order of vanishing modulo $\ell$.	In \cite{H} we studied
the reduction when $\E_{f,\ell}$ has `small' monodromy instead.

The usual Weil pairing on $E_f[\ell]\times E_f[\ell]$ extends to a
non-degenerate alternating pairing $\E_{f,\ell}\times\E_{f,\ell}\to\Fell(1)$,
hence $\E_{f,\ell}\to C$ is self-dual.  Together with Poincar\'e duality we
obtain a non-degenerate {\it symmetric} pairing of cohomology groups
    $$ H^1(C\times\Fqbar,\E_{f,\ell}) \times H^1(C\times\Fqbar,\E_{f,\ell})
	\longrightarrow H^2(C\times\Fqbar,\Fell(1)).
    $$
That is, we have a non-degenerate orthogonal pairing of $\Fell$-vector
spaces $V_{f,\ell} \times V_{f,\ell} \to \Fell.$  We write $O(V_{f,\ell})$
for the subgroup of $\GL(V_{f,\ell})$ preserving the pairing and note
that $\Frob_q$ preserves the pairing on $V_{f,\ell}$, so belongs to a
well-defined conjugacy class in $O(V_{f,\ell})$.

\subsection{Families of Twists}

For ease of exposition we assume $C=\P^1$ and $K=\F_q(t)$.  For general
$C$ one will have to increase the minimum value of $\ell$ required for
a surjectivity statement of the form of theorem~\ref{thm1}.

We fix an elliptic curve $E_1/K$ such that $\E_1\to C$ has at least
one fiber of multiplicative reduction away from $\infty$.  We also
fix a non-zero polynomial $m\in\F_q[t]$ which vanishes at one or more
finite points in $M_1$ in order that for every $f\in\F_q[t]$ which is
relatively prime to $m$, the twist $\E_f\to C$ also has at least one
fiber of multiplicative reduction away from $\infty$.

For every integer positive integer $d$ we define the family of twisting
polynomials
    $$
	F_d = \left\{
	    f\in\Fqbar[t] :
	    \mbox{$f$ is square-free}, \deg(f)=d, \gcd(f,m) = 1
	\right\}.
    $$
We may regard $F_d$ as a $(d+1)$-dimensional affine space and we write
$F_d(\F_{q^n})$ for the set of $\F_{q^n}$-valued points, i.e. the subset
of $f$ with coefficients in $\F_{q^n}$.  Unless we restrict the leading
coefficient of $f$, there will be many $g\in F_d$ which give rise to an
isomorphic twist.  However, for every $n\geq 1$, the number of twists
$\E_g$ by $g\in F_d(\F_{q^n})$ isomorphic to $\E_f$ is independent of
$f\in F_d(\F_{q^n})$.

Katz first suggested restricting to twists parametrized by a fixed $F_d$
in part because the set of twists by $f\in F_d(\F_{q^n})$ satisfy
a remarkable uniformity property: the degree of the $L$-function
$L(T,E_f/K_n)$ is independent of $f$ and $n$.  In fact, this is
a consequence of a deeper sheaf-theoretic statement: for every
non-exceptional $\ell$, there is a unique \etale $\Fell$-lisse sheaf
$\T_{d,\ell}\to F_d$ whose (geometric) fiber over any $f\in F_d(\F_{q^n})$
is the $\Fell$-vector space $H^1(C\times\overline{\F}_{q^n},\E_{f,\ell})$.

We fix a geometric point $f\in F_d$ and let $\pi_1(F_d)=\pi_1(F_d,f)$
denote the \etale fundamental group.  Then $\T_{d,\ell}$ corresponds
to a $\Fell$-representation $\rho:\pi_1(F_d)\to\GL(V_{f,\ell})$ and
the image is well defined up to inner automorphism.  We define the
arithmetic monodromy group to be the image of $\pi_1(F_d)$ and the
geometric monodromy group to be the image of $\pi_1(F_d\times\Fqbar)$.
If we take $\T_{d,\ell}$ together with its orthogonal pairing, then the
results at the end of the previous section imply that both monodromy
groups lie in $O(V_{f,\ell})$.

The main question of interest for us is to determine as precisely as
possible the monodromy groups as $d$ and $\ell$ vary.  While we do not
answer this question completely, the following theorem demonstrates that
the monodromy is usually `big' in a strongly uniform way.

\begin{theorem}\label{thm1}
    If $\ell\geq 5$ and $\deg(f)\gg_E 0$, then the geometric monodromy group
    has index at most two in $O(V_{f,\ell})$ and is not $SO(V_{f,\ell})$.
    That is, $G$ is one of the following:
    \begin{enumerate}
    \item the full orthogonal group $O(V_{f,\ell})$;
    \item the kernel of the spinor norm;
    \item the kernel of the product of the spinor norm and the determinant.
    \end{enumerate}
\end{theorem}

In order to prove the theorem it suffices to restrict to one-parameter
families parametrized by an affine curve $U\subset F_d$ and to show that
the image of $\pi_1(U\times\Fqbar)$ is already big.  More precisely, fix
any $g\in F_{d-1}$ and consider the one-parameter family of polynomials
$(c-t)\cdot g(t)$.  We let $U_g\subset\P^1$ be the open subset of
$c\in\A^1$ such that $(c-t)\cdot g(t)\in F_d$.  Then theorem~\ref{thm1}
follows immediately from the following theorem whose proof we postpone
until the next section.

\begin{theorem}\label{thm2}
    If $\ell\geq 5$ and $\deg(g)\gg_E 0$, then the image of
    $\pi_1(U_g\times\Fqbar)$ has index at most two in $O(V_{f,\ell})$ and
    is not $SO(V_{f,\ell})$.
\end{theorem}

\subsection{Katz One-Parameter Families of Twists}\label{sec::katz}

We fix $g\in F_{d-1}$ and let $U_g\subset F_d$ be as before.
The key observation Katz makes to prove analogous $\ell$-adic
monodromy theorems is that the restriction of $\T_{d,\ell}$ to $U_g$
is the middle convolution sheaf $\MC_{-1}(\E_g)\to C$ (cf.~section~\ref{sec::conv}).
In particular, $\T_{d,\ell}$ is irreducible and tame and we can describe
its monodromy around the points of $\P^1-U_g$.  We refer the reader
to section~\ref{sec::group} for the definition of a reflection and
isotropic shear.

\begin{lemma}\label{lemma10}
    For every geometric point $c\in\A^1-U_g$ fix a
    topological generator $\sigma_c$ of the inertia group
    $I(c)\leq\pi_1^{\tame}(U_g\times\Fqbar)$ and let $V=V_{f,\ell}$.
    \begin{enumerate}
    \item If $\E_g\to C$ has good reduction over $t=c$, then $\sigma_c$
	acts trivially on $V$.
    \item If $\E_g\to C$ has multiplicative reduction over $t=c$, then
	$\sigma_c$ acts as a reflection on $V$.
    \item If $\E_g\to C$ has additive reduction of Kodaira type $I_0^*$
	over $t=c$, then $\sigma_c$ acts as an isotropic shear on $V$.
    \end{enumerate}
    For all other $c\in\A^1-U_g$, $\sigma_c$ acts as a non-scalar on
    the two-dimensional quotient $V/V^{\sigma_c=1}$.
\end{lemma}

\begin{proof}
    We follow the notation of section~\ref{sec::conv}.  The fiber of the
    convolution $\MC_{-1}(\E_g)$ over a (geometric) closed point $c\in
    U_g$ is $H^1(C,i_*i^*(\E_g\otimes\tau_c^*\L_{-1}))$, where $i:U_g\to
    C$.  One can easily show that $i_*i^*(\E_g\otimes\tau_c^*\L_{-1})$
    is the fiber of $\E_f$ over $c$, hence the restrictions
    $\T_{d,\ell}\to U_g$ and $\MC_{-1}(\E_g)\to U_g$ are isomorphic,
    so $\T_{d,\ell}\simeq\MC_{-1}(\E_g)$.

    If $\E_g\to C$ has multiplicative reduction over $t=c$,
    then $\E_g(c)/\E_g(c)^{I(c)}$ is the trivial representation
    $\F_\ell$, so $\T_{d,\ell}(c)/\T_{d,\ell}(c)^{I(c)}$
    is the scalar representation $-1$ and the monodromy is a
    reflection.  If $\E_g\to C$ has additive reduction of Kodaira
    type $I_0^*$ over $t=c$, then $\E_g(c)/\E_g(c)^{I(c)}$ is the
    two-dimensional representation $(\F_\ell\oplus\F_\ell)\otimes -1$,
    so $\T_{d,\ell}(c)/\T_{d,\ell}(c)^{I(c)}$ is $\F_\ell\oplus\F_\ell$.
    Thus $\T_{d,\ell}$ has two unipotent blocks of the form $U_2$ and all
    other blocks are trivial (cf.~(\ref{eqn1})).  Finally, for all other
    types of (additive) reduction we see that $\E_g(c)/\E_g(c)^{I(c)}$
    is two-dimensional and $\sigma_c$ acts as a non-scalar, so the same
    is true for $\T_{d,\ell}(c)/\T_{d,\ell}(c)^{I(c)}$.
\end{proof}

The elements $S=\{\sigma_c:c\in\A^1-U_g\}$ topologically generate
$\pi_1^{\tame}(U_g\times\Fqbar)$, hence they generate the image $G$
of $\pi_1^{\tame}(U_g\times\Fqbar)$ in $O(V)$.  One implication is that
we can ignore the monodromy around $\infty$ (which is more difficult to
describe).  More importantly, the codimension of $V^{\sigma=1}$ in $V$ is
at most 2 for every $\sigma\in S$, and this places severe restrictions on
the monodromy when $|S|\gg 0$.  In particular, as $\deg(g)$, hence $|S|$,
tends to infinity, the complement of the subset of isotropic shears in
$S$ has bounded order.  Therefore theorem~\ref{thm2} is a consequence
of theorem~\ref{thm3}, applied with $r=2$ and $S_0$ the complement of
the reflections and the elements with order prime to $(r+1)!=6$.

\subsection{Example: Twists of the Legendre Curve}

Let $K=k(\lambda)$ and let $E/K$ be the elliptic curve with affine
model $y^2=x(x-1)(x-\lambda)$.  We write $F_d$ for the square-free
polynomials in $\kbar[\lambda]$ of degree $d$ which are relatively prime
to $\lambda(\lambda-1)$ and $F_d(k)$ for the subset in $k[\lambda]$.  For
each $f\in F_d(k)$ we write $M_f,A_f$ for the divisors of multiplicative,
additive reduction of the N\'eron model $\E_f\to C$ and $\div_0(f)$ for
the divisor of zeros of $f$.
\begin{lemma} Suppose $f\in F_d(k)$.  Then
$$ M_f,A_f,\dim(V_{f,\ell})\ =\ \begin{cases}
    \,\ \ \{0,1\},\ \ \div_0(f)\cup\{\infty\},\ \ 2d\ \ \ 
	& \mbox{if $d$ is even} \\
    \,\{0,1,\infty\},\qquad\div_0(f),\quad 2d-1\,  & \mbox{if $d$ is odd}
\end{cases}. $$
Moreover, the support of the fibers of $\E_f\to C$ of Kodaira type $I_0^*$
is $\div_0(f)$.
\end{lemma}

\begin{proof}
    Everything except the dimension assertions are proved in lemma 7 of
    \cite{H}, while $\dim(V_{f,\ell})$ is given in lemma~\ref{lemma11}.
\end{proof}

By theorem~\ref{thm4} the sheaves $\E_{f,\ell}\to C$ have big monodromy
for all $f\in F_d(k)$ and all odd $\ell$, so for $g\in F_{d-1}(k)$,
$\ell$ odd, and $U_g=C-M_g-A_g$, the restriction $\T_{d,\ell}\to U_g$
is irreducible and lisse.  The monodromy about a geometric point in
$M_g-\{\infty\}$ is a reflection and the monodromy about a geometric
point in $A_g-\{\infty\}$ is an isotropic shear by lemma~\ref{lemma10},
so $S_0=\emptyset$ (if $\ell>3$) and theorem~\ref{thm1} holds for $d\geq
2$ and $\ell\geq 5$.

One can derive similar results if we replace $E/K$ by the twist
$E_\lambda/K$.  In particular, we can construct $V_{f,\ell}$ satisfying
$\dim(V_{f,\ell})=2d+1\equiv 3\pmod{4}$ for $d$ odd.  In order to construct
examples with $\dim(V_{f,\ell})\equiv 2\pmod{4}$ one should replace the
Legendre curve by one of the curves denoted $X_{211},X_{321},X_{431}$ in
\cite{MP}.  In particular, up to an automorphism of the base $C=\P^1$, we
can assume $M_1=\{0,1\}$ and $A_1=\{\infty\}$ as before, but now the key
difference is that $\infty\in A_f$ for every twist.

\rmk Together these examples increase the set of big subgroups of
orthogonal groups which are known to occur as Galois groups over
$\Q(t)$.  It is difficult to say precisely which group occurs, but
despite the ambiguity these extend previous results (cf.~survey in
\cite{MM}).

\subsubsection{Almost Independent}

In this section we assume $k$ is finite or separably closed.  Fix
$g\in F_{d-1}(k)$ and let $\T_{d,\ell}\to U_g$ be as in the previous
section.  Let $L$ denote the function field $k(U_g)$ and let $L_\ell$
denote the splitting field $L(V_{f,\ell})$.  As the following theorem
shows, up to replacing $L$ by a finite extension, these extensions
are almost independent (cf.~10.1?~of \cite{S4}).

\begin{theorem}
    If $d\geq 3$, then there is a finite extension $M/L$ so
    $L_{\ell_1}\cap L_{\ell_2}\leq M$ for $\ell_1>\ell_2\geq 5$.
\end{theorem}

\begin{proof}
    By the results in the previous section, if $d\geq 2$ and $\ell\geq 5$,
    then $G_\ell=\Gal(L_\ell/L)$ is a big subgroup of an orthogonal group
    $\Gamma_\ell$ and $Q_\ell=G_\ell/\D G_\ell$ is a subgroup of
    $\Z/2\oplus\Z/2$.  Moreover, if $\ell_1>\ell_2\geq 5$ and $d\geq 3$ (so
    $\dim(V_{f,\ell})\geq 5$), then the quotients $\D G_{\ell_i}/Z(\D
    G_{\ell_i})$ are non-abelian, simple, and pairwise non-isomorphic for
    $i=1,2$ (cf.~theorem 5.27 of \cite{Ar}).  Therefore Ribet's lemma
    (5.2.2 of \cite{R}) implies $L_{\ell_1}\cap L_{\ell_2}$ is contained in
    the fixed fields of $Q_{\ell_1},Q_{\ell_2}$.  On the other hand, the
    fixed field of $Q_\ell$ corresponds to an unramified cover
    $V_{g,\ell}\to U_g$ of bounded degree.  In particular, there are only
    finitely many covers which occur as we vary $\ell$ because of our
    assumptions on $k$, so we make take $M$ to be the compositum of all the
    corresponding extensions.
\end{proof}

\subsection{Generalizing to Abelian Varieties}

While the previous sections deal exclusively with twists of elliptic
curves, both for ease of exposition and application, most of the results
can be easily adapted to deal with twists of `many' abelian varieties
$A_1/K$ of dimension $g$ with trivial $K/k$-trace.  The easiest is
to assume that the Galois group $G_\ell$ of $K(A_1[\ell])/K$ is big for
$\ell\gg_{A_1} 0$, but it suffices to assume it acts irreducibly on $A_1[\ell]$
for $\ell\gg_{A_1} 0$.  Either way we must also assume that $G_\ell$ acts
tamely on $A_1[\ell]$; this is automatic if the characteristic of $K$
is sufficiently large with respect to the genus of $C$.  We must also
assume that $A_1/K$ has at potentially semi-stable reduction with toric
part of dimension one over some closed point $x\in C$.  For example,
we may take $A_1/K$ to be any $J/K$ from section~\ref{sec::yu}.

Over almost all the closed points in $C$ an arbitrary quadratic twist of
$A_1/K$ has either good reduction or totally additive reduction.  In the
latter case \cite{LO} implies the component group of the fiber over $x$
of the \Neron model is uniformly bounded by a constant which depends
only on the dimension of $A_1$.  For each of the remaining closed points
the component group of the special fiber belongs to a set of at most two
finite groups.  Therefore as we vary over the quadratic twists $A_f/K$
of $A_1/K$, the set of primes dividing the order of the component group
of the \Neron model $\AA_f\to C$ is finite.  In particular, if $\Lambda$
is sufficiently small, then for every $\ell\in\Lambda$, we may assume
$G_\ell$ acts irreducibly on $A_1[\ell]$ and $\ell$ is relatively prime
to the order of the component group of $\AA_f\to C$ for every twist.

We make these assumptions because they imply the cohomology groups
of the $\Zell$-sheaf $T_\ell(\AA_f)\to C$, the latter defined as the
projective system of the \etale sheaves $\AA_{f,\ell^n}$ (cf.~section
2.2 of \cite{G2}), are sufficiently well behaved.  In particular, the
$\Zell$-sheaf $\T_{d,\ell^\infty}\to C$, defined as the projective system
of the generalized sheaves $\T_{d,\ell^n}\to C$, is torsion free and
$\T_{d,\ell^\infty}\otimes_{\Zell}\Fell$ is isomorphic to $\T_{d,\ell}$.
The key is to consider the Kummer sequence
    $$ 0\longrightarrow T_\ell(\AA_f)
        \overset{\times\ell}\longrightarrow T_\ell(\AA_f)
        \longrightarrow \AA_{f,\ell}
        \longrightarrow 0 $$
which is defined as the projective system of the sequences
    $$ 0\longrightarrow \AA_{f,\ell^n}
        \longrightarrow \AA_{f,\ell^{n+1}}
        \overset{\times\ell^n}\longrightarrow \AA_{f,\ell}
        \longrightarrow 0,
    \qquad n\geq 0. $$
The corresponding cohomology sequence simplifies (cf.~2.1--2.4 of
\cite{Shi}) to
$$   0  \longrightarrow H^1(\Cbar,T_\ell(\AA_f))
		\overset{\times\ell}\longrightarrow H^1(\Cbar,T_\ell(\AA_f))
        \longrightarrow H^1(\Cbar, \AA_{f,\ell})
        \longrightarrow 0 $$
because $H^i(\Cbar,\AA_{f,\ell})=0$ for $i\neq 1$, which implies the
claim.

\end{document}